\documentclass{amsart}
\usepackage{hyperref}

\newtheorem{theorem}{Theorem}[section]
\newtheorem{lemma}[theorem]{Lemma}
\newtheorem{corollary}[theorem]{Corollary}

\newcommand{\R}{{\mathbb R}}
\newcommand{\N}{{\mathbb N}}
\newcommand{\Z}{{\mathbb Z}}
\newcommand{\C}{{\mathbb C}}

\newcommand{\nn}{\nonumber}
\newcommand{\be}{\begin{equation}}
\newcommand{\ee}{\end{equation}}
\newcommand{\bea}{\begin{eqnarray}}
\newcommand{\eea}{\end{eqnarray}}

\newcommand{\ol}{\overline}

\newcommand{\spr}[2]{\langle #1 , #2 \rangle}
\newcommand{\id}{{\rm 1\hspace{-0.6ex}l}}

\newcommand{\lz}{\ell^2(\Z)}
\newcommand{\tl}{\mathrm{TL}}
\newcommand{\tr}{\mathrm{tr}}

\newcommand{\Rg}[1]{R_{2g+2}^{1/2}(#1)}

\newcommand{\vprod}[2]{\!\!\!\!\begin{array}{c} \mbox{\raisebox{-0.5ex}[0.5ex]
{$\scriptstyle #2 $}} \\ \displaystyle \hspace*{1.1ex}\prod{}^* \\
\mbox{\raisebox{0.6ex}[-0.6ex]{$ \scriptstyle  #1 $}} \end{array}}
\newcommand{\vsum}[2]{\!\!\!\!\begin{array}{c} \mbox{\raisebox{-0.5ex}[0.5ex]
{$\scriptstyle #2 $}} \\ \displaystyle \hspace*{1.1ex}\sum{}^* \\
\mbox{\raisebox{0.6ex}[-0.6ex]{$ \scriptstyle  #1 $}} \end{array}}

\newcommand{\sig}{\sigma}

\numberwithin{equation}{section}

\begin{document}

\title[Trace Formulas in Connection with Scattering Theory]{Trace Formulas in Connection with Scattering Theory for Quasi-Periodic Background}

\author{Johanna Michor}
\address{Fakult\"at f\"ur Mathematik\\
Nordbergstrasse 15\\ 1090 Wien\\ Austria\\ and International Erwin Schr\"odinger
Institute for Mathematical Physics, Boltzmanngasse 9\\ 1090 Wien\\ Austria}
\email{\href{mailto:Johanna.Michor@esi.ac.at}{Johanna.Michor@esi.ac.at}}
\urladdr{\href{http://www.mat.univie.ac.at/~jmichor/}{http://www.mat.univie.ac.at/\~{}jmichor/}}

\author{Gerald Teschl}
\address{Institut f\"ur Mathematik\\
Nordbergstrasse 15\\ 1090 Wien\\ Austria\\ and International Erwin Schr\"odinger
Institute for Mathematical Physics, Boltzmanngasse 9\\ 1090 Wien\\ Austria}
\email{\href{mailto:Gerald.Teschl@univie.ac.at}{Gerald.Teschl@univie.ac.at}}
\urladdr{\href{http://www.mat.univie.ac.at/~gerald/}{http://www.mat.univie.ac.at/\~{}gerald/}}

\thanks{Work supported by the Austrian Science Fund (FWF) under Grant
No.\ P17762}


\begin{abstract}
We investigate trace formulas for Jacobi operators which are trace class 
perturbations of quasi-periodic finite-gap operators using Krein's
spectral shift theory. In particular we establish the conserved quantities for the
solutions of the Toda hierarchy in this class.
\end{abstract}

\maketitle

\section{Introduction}

Scattering theory for Jacobi operators $H$ with periodic (respectively more general)
background has attracted considerable interest recently. In \cite{voyu} Volberg and
Yuditskii have treated the case where $H$ has a homogeneous spectrum and is of
Szeg\"o class exhaustively. In \cite{emtqps} Egorova and we have established direct
and inverse scattering theory for Jacobi operators which are short range perturbations
of quasi-periodic finite-gap operators. For further information and references we refer
to these articles and \cite{tjac}.

In the case of constant background it is well-known that the transmission coefficient is
the perturbation determinant in the sense of Krein \cite{krein}, see e.g., \cite{tist} or
\cite{tjac}. The purpose of the present paper is to establish this result  for the
case of quasi-periodic finite-gap background, thereby establishing the connection with
Krein's spectral shift theory. 

Moreover, scattering theory for Jacobi operators is not only interesting in its own
right, it also constitutes the main ingredient of the inverse scattering transform for
the Toda hierarchy (see, e.g., \cite{fl2}, \cite{fad}, \cite{tjac}, or \cite{ta}).
Since the transmission coefficient is invariant when our Jacobi
operator evolves in time with respect to some equation of the Toda hierarchy, the
corresponding trace formulas provide the conserved quantities for the
Toda hierarchy in this setting.

\section{Notation}

Let
\be
H_q f(n) = a_q(n) f(n+1) + a_q(n-1) f(n-1) + b_q(n) f(n)
\ee
be a quasi-periodic Jacobi operator in $\lz$ associated with the Riemann surface
of the function
\be
\Rg{z}, \qquad R_{2g+2}(z) = \prod_{j=0}^{2g+1} (z-E_j), \qquad
E_0 < E_1 < \cdots < E_{2g+1},
\ee
$g\in \N$. The spectrum of $H_q$ is purely absolutely continuous and
consists of $g+1$ bands
\be
\sig(H_q) = \bigcup_{j=0}^g [E_{2j},E_{2j+1}].
\ee
For every $z\in\C$ the Baker-Akhiezer functions $\psi_{q,\pm}(z,n)$ are
two (weak) solutions of $H_q \psi = z \psi$, which are linearly independent away
from the band-edges $\{E_j\}_{j=0}^{2g+1}$, since their Wronskian is given by
\be
W_q(\psi_{q,-}(z), \psi_{q,+}(z)) = \frac{R^{1/2}_{2g+2}(z)}{\prod_{j=1}^g (z-\mu_j)}.
\ee
Here $\mu_j$ are the Dirichlet eigenvalues at base point $n_0=0$.
We recall that $\psi_{q,\pm}(z,n)$ have the form
\[
\psi_{q,\pm}(z,n) =\theta_{q,\pm}(z,n) w(z)^{\pm n},
\]
where $\theta_{q,\pm}(z,n)$ is quasi-periodic with respect to $n$ and
$w(z)$ is the quasi-momentum. In particular, $|w(z)|<1$ for
$z\in\C\backslash\sig(H_q)$ and $|w(z)|=1$ for $z\in\sig(H_q)$.

We assume that the reader is familiar with this class of operators and refer to
\cite{emtqps} and \cite{tjac} for further information.

\section{Asymptotics of Jost solutions}
\label{secJS}

After we have these preparations out of our way, we come to
the study of short-range perturbations $H$ of $H_q$ associated with sequences $a$, $b$
satisfying $a(n) \rightarrow a_q(n)$ and $b(n) \rightarrow b_q(n)$ as $|n|
\rightarrow \infty$.  More precisely, we will make the following assumption throughout
this paper:

Let $H$ be a perturbation of $H_q$ such that
\be                         \label{hypo}
\sum_{n \in \mathbb{Z}} \Big(|a(n) - a_q(n)| + |b(n) - b_q(n)| \Big) <
\infty,
\ee
that is, $H-H_q$ is trace class.

We first establish existence of Jost solutions, that is, solutions of the
perturbed operator which asymptotically look like the Baker-Akhiezer solutions.

\begin{theorem} \label{thmjost}
Assume (\ref{hypo}). Then there exist (weak) solutions
$\psi_{\pm}(z, .)$, $z \in \C\backslash\{E_j\}_{j=0}^{2g+1}$, of
$H \psi = z \psi$ satisfying
\be \label{jost1}
\lim_{n \rightarrow \pm \infty}
w(z)^{\mp n} \left( \psi_{\pm}(z, n) - \psi_{q,\pm}(z, n) \right) = 0,
\ee
where $\psi_{q,\pm}(z, .)$ are the Baker-Akhiezer functions.
Moreover, $\psi_{\pm}(z, .)$ are continuous 
(resp.\ holomorphic) with respect to $z$ whenever $\psi_{q,\pm}(z, .)$ are and
have the following asymptotic behavior
\be                 \label{B4jost}
\psi_\pm(z,n) =   \frac{z^{\mp n}}{A_\pm(n)} \Big(\vprod{j=0}{n-1}a_q(j)\Big)^{\pm 1} 
\Big(1 + \Big(B_\pm(n) \pm \vsum{j=1}{n} b_q(j- {\scriptstyle{0 \atop 1}}) \Big)\frac{1}{z}
+ O(\frac{1}{z^2}) \Big),
\ee
where
\bea \nn
A_+(n) &=& \prod_{j=n}^{\infty} \frac{a(j)}{a_q(j)}, \qquad
B_+(n)= \sum_{m=n+1}^\infty (b_q(m)-b(m)), \\
A_-(n) &=& \prod_{j=- \infty}^{n-1} \frac{a(j)}{a_q(j)}, \qquad
B_-(n) = \sum_{m=-\infty}^{n-1} (b_q(m)-b(m)).
\eea
\end{theorem}

\begin{proof}
The proof can be done as in the periodic case (see e.g., \cite{emtqps}, \cite{gerass},
\cite{tosc} or \cite{tjac}, Section 7.5). There a stronger decay assumption (i.e.,
first moments summable) is made, which is however only needed at the band edges
$\{E_j\}_{j=0}^{2g+1}$.
\end{proof}

\noindent
For later use we note the following immediate consequence

\begin{corollary} \label{corpisprime}
Under the assumptions of the previous theorem we have
\be \label{jost2}
\lim_{n \rightarrow \pm \infty}
w(z)^{\mp n} \left( \psi_{\pm}'(z, n) \mp n \frac{w'(z)}{w(z)} \psi_{\pm}(z, n) -
\psi_{q,\pm}'(z, n) \pm n \frac{w'(z)}{w(z)} \psi_{q,\pm}(z, n) \right) = 0,
\ee
where the prime denotes differentiation with respect to $z$.
\end{corollary}

\begin{proof}
Just differentiate (\ref{jost1}) with respect to $z$, which is permissible by uniform
convergence on compact subsets of $\C\backslash \{E_j\}_{j=0}^{2g+1}$.
\end{proof}

\noindent
We remark that if we require our perturbation to satisfy the usual short range
assumption as in \cite{emtqps} (i.e., first moments summable), then we even have
$w(z)^{\mp n} (\psi_{\pm}'(z, n) - \psi_{q,\pm}'(z, n)) \to 0$.

 From Theorem~\ref{thmjost} we obtain a complete characterization of the spectrum of $H$.

\begin{theorem} 
Assume (\ref{hypo}). Then we have $\sig_{ess}(H)=\sig(H_q)$, the
point spectrum of $H$ is confined to $\ol{\R\backslash\sig(H_q)}$. Furthermore,
the essential spectrum of $H$ is purely absolutely continuous except for possible
eigenvalues at the band edges.
\end{theorem}

\begin{proof}
An immediate consequence of the fact that $H-H_q$ is trace class and boundedness
of the Jost solutions inside the essential spectrum.
\end{proof}

\noindent
Our next result concerns the asymptotics of the Jost solutions at the {\em other side}.

\begin{lemma} \label{lemothers}
Assume (\ref{hypo}). Then the Jost solutions
$\psi_{\pm}(z, .)$, $z \in \C\backslash\sig(H)$, satisfy
\be                         \label{perturbed sol}
\lim_{n \rightarrow \mp \infty}
|w(z)^{\mp n} (\psi_{\pm}(z, n) - \alpha(z)\psi_{q,\pm}(z, n))| = 0,
\ee
where
\be                  
\alpha(z) = \frac{W(\psi_-(z),\psi_+(z))}{W_q(\psi_{q,-}(\lambda), \psi_{q,+}(z))} =
\frac {\prod_{j=1}^g(z - \mu_j)}{\Rg{z}} W(\psi_-(z), \psi_+(z)).
\ee
\end{lemma}

\begin{proof}
Since $H-H_q$ is trace class, we have for the difference of the Green's functions
\[
\lim_{n\to\pm\infty} G(z,n,n)- G_q(z,n,n) =
\lim_{n\to\pm\infty} \spr{\delta_n}{((H-z)^{-1} - (H_q-z)^{-1}) \delta_n} =0
\]
and hence
\[
\lim_{n\to-\infty} \psi_{q,-}(z,n)(\psi_+(z,n) - \alpha(z) \psi_{q,+}(z,n)) =0,
\]
which is the claimed result.
\end{proof}

\noindent
Note that $\alpha(z)$ is just the inverse of the transmission coefficient (see, e.g.,
\cite{emtqps} or \cite{tjac}, Section~7.5).
It is holomorphic in $\C\backslash\sig(H_q)$ with simple zeros at the discrete eigenvalues
of $H$ and has the following asymptotic behavior
\be
\alpha(z) = \frac{1}{A} ( 1 + \frac{B}{z} + O(z^{-2})), \qquad
A=A_-(0)A_+(0), \quad B= B_-(1)+B_+(0).
\ee

\section{Connections with Krein's spectral shift theory and Trace formulas}

To establish the connection with Krein's spectral shift theory we next
show:

\begin{lemma}
We have
\be
\frac{d}{dz} \alpha(z) = - \alpha(z) \sum_{n \in \Z} \big( G(z, n, n) - G_q(z, n, n)\big),
\qquad z\in\C\backslash\sig(H),
\ee
where $G(z,m,n)$ and $G_q(z,m,n)$ are the Green's functions of $H$ and $H_q$,
respectively.
\end{lemma}

\begin{proof}
Green's formula (\cite{tjac}, eq. (2.29)) implies
\be     \label{green 1}
W_n(\psi_+(z), \psi_-'(z)) - W_{m-1}(\psi_+(z), \psi_-'(z)) =  
\sum_{j=m}^n \psi_+(z,j) \psi_-(z,j), 
\ee 
hence the derivative of the Wronskian can be written as
\bea \nn
\lefteqn{\frac{d}{dz}W(\psi_-(z), \psi_+(z)) = W_n(\psi_-'(z), \psi_+(z)) +
W_n(\psi_-(z), \psi_+'(z))} \\ \nn
&=& W_m(\psi_-'(z), \psi_+(z)) + W_n(\psi_-(z), \psi_+'(z)) -
\sum_{j=m+1}^n \psi_+(z,j)\psi_-(z,j).
\eea
Using Corollary~\ref{corpisprime} and Lemma~\ref{lemothers} 
we have 
\bea \nn
W_m(\psi_-'(z), \psi_+(z)) &=& W_m(\psi_-' + m \frac{w'}{w} \psi_-, \psi_+) -\\
\nn
&& \frac{w'}{w} \big( m\, W(\psi_-, \psi_+) - a(m) \psi_-(m+1) \psi_+(m) \big)\\ \nn
&\to& \alpha W_{q,m} (\psi_{q,-}' + m \frac{w'}{w} \psi_{q,-}, \psi_{q,+}) -\\ \nn
&& \alpha \frac{w'}{w} \big( m\, W_q(\psi_{q,-}, \psi_{q,+}) -
a_q(m) \psi_{q,-}(m+1) \psi_{q,+}(m) \big)\\ \nn
&=& \alpha(z) W_m(\psi_{q,-}'(z), \psi_{q,+}(z))
\eea
as $m \rightarrow - \infty$. Similarly we obtain
\bea \nn
W_n(\psi_-(z), \psi_+'(z)) &\to& \alpha(z)
W_n(\psi_{q,-}(z), \psi_{q,+}^{\prime}(z))
\eea
as $n \rightarrow \infty$
and again using (\ref{green 1}) we have
\[
W_m(\psi_{q,-}^{\prime}(z), \psi_{q,+}(z)) 
= W_n(\psi_{q,-}^{\prime}(z), \psi_{q,+}(z))
+ \sum_{j=m+1}^n \psi_{q,+}(z,j) \psi_{q,-}(z,j).
\] 
Collecting terms we arrive at
\bea \nn
W^{\prime}(\psi_-(z), \psi_+(z)) &=& 
- \sum_{j \in \Z} \Big( \psi_+(z, j) \psi_-(z, j) - 
\alpha(z) \psi_{q,+}(z, j) \psi_{q,-}(z, j) \Big) \\  \nn
&& + \alpha(z) W_q^{\prime}(\psi_{q,-}(z) \psi_{q,+}(z)).
\eea
Now we compute
\bea \nn
\frac{d}{dz} \alpha(z) &=& \frac{d}{dz} \Big( \frac{W}{W_q}\Big) = 
\Big(\frac{1}{W_q}\Big)^{\prime} W + \frac{1}{W_q} W^{\prime} \\ \nn 
&=& - \frac{W_q^{\prime}}{W_q^2} W + \frac{1}{W_q}
\Big( - \sum_{j \in \Z} \Big( \psi_+  \psi_-  - \alpha  \psi_{q,+}  \psi_{q,-} \Big)    
+ \alpha  W_q^{\prime}\Big) \\ \nn
&=& - \frac{1}{W_q} \sum_{j \in \Z} \Big( \psi_+  \psi_-  - \alpha  \psi_{q,+}  \psi_{q,-} \Big), 
\eea
which finishes the proof.
\end{proof}

\noindent
As an immediate consequence, we can identify $\alpha(z)$ as Krein's perturbation determinant (\cite{krein}) of the pair $H$, $H_q$.

\begin{theorem} 
The function $A \alpha(z)$ is Krein's perturbation determinant:
\be
\alpha(z) = \frac{1}{A} \det \big(\id + (H(t)-H_q(t)) (H_q(t)-z)^{-1}\big), \quad
A = \prod_{j=- \infty}^\infty \frac{a(j)}{a_q(j)}.
\ee
\end{theorem}

\noindent
By \cite{krein}, Theorem~1, $\alpha(z)$ has 
the following representation 
\be
\alpha(z) = \frac{1}{A} \exp \Big(\int_{\R} 
\frac{\xi_{\alpha}(\lambda)d\lambda}{\lambda - z} \Big),
\ee
where
\be
\xi_{\alpha}(\lambda) = \frac{1}{\pi}\lim_{\epsilon \downarrow 0} 
\arg \alpha(\lambda + i \epsilon) 
\ee
is the spectral shift function.

Hence
\be
\tau_j = \tr(H^j - (H_q)^j) = j \int_{\R} \lambda^{j-1} \xi_{\alpha}(\lambda)d\lambda,
\ee
where $\tau_j/j$ are the expansion coefficients of $\ln \alpha(z)$ around $z=\infty$:
\[
\ln \alpha(z) = -\ln A  - \sum_{j=1}^\infty \frac{\tau_j}{j\, z^j}.
\]
They are related to the expansion $\alpha_j$ coefficients of
\[
\alpha(z) = \frac{1}{A} \sum_{j=0}^\infty \frac{\alpha_j}{z^j}, \qquad \alpha_0=1,
\]
via
\be
\tau_1 = - \alpha_1, \qquad \tau_j = - j \alpha_j - \sum_{k=1}^{j-1} \alpha_{j-k} \tau_k.
\ee

\section{Conserved quantities of the  Toda hierarchy}

Finally we turn to solutions of the Toda hierarchy $\tl_r$ (see, e.g., \cite{bght}, \cite{fad},
\cite{tjac}, or \cite{ta}). Let $(a_q(t),b_q(t))$ be a quasi-periodic finite-gap solution of
some equation in the Toda hierarchy, $\tl_r(a_q(t),b_q(t))=0$, and let $(a(t),b(t))$
be another solution, $\tl_r(a(t),b(t))=0$, such that (\ref{hypo}) holds for all $t$.

Since the transmission coefficient $T(z,t)=T(z,0)\equiv T(z)$ is conserved (see
\cite{emtist} -- formally this follows from unitary invariance of the determinant),
so is $\alpha(z)= T(z)^{-1}$.

\begin{theorem}
The quantities
\be
A = \prod_{j=- \infty}^\infty \frac{a(j,t)}{a_q(j,t)}
\ee
and $\tau_j= \tr(H^j(t) - H_q(t)^j)$, that is,
\bea \nn
\tau_1 &=& \sum_{n\in\Z} b(n,t) - b_q(n,t)\\ \nn
\tau_2 &=& \sum_{n\in\Z} 2( a(n,t)^2 - a_q(n,t)^2) + (b(n,t)^2 - b_q(n,t)^2)\\ \nn
&\vdots&
\eea
are conserved quantities for the Toda hierarchy.
\end{theorem}

\end{document}